21 February, 1995

# The Artinian Berger Conjecture

Guillermo Cortiñas  
Dept. de Matemática  
Universidad de La Plata  
50 y 115  
(1900) La Plata  
Argentina

Susan C. Geller[1]  
Dept. of Mathematics  
Texas A&M University  
College Station, TX 77802  
USA

Charles A. Weibel[2]  
Dept. of Mathematics  
Rutgers University  
New Brunswick, NJ 08903  
USA

**Abstract.** We propose an Artinian version of Berger's Conjecture for curves, concerning the module of Kähler differentials of an algebra. Our version implies Berger's Conjecture in characteristic 0. We establish our Artinian Berger Conjecture in a number of cases, and prove that Berger's Conjecture holds for curve singularities whose conductor ideal contains the cube of a maximal ideal.

In this paper we introduce and study a conjecture, which we call the Artinian Berger Conjecture (or "ABC"[3]), about the Kähler differentials $\Omega_{A/k}$ of a finite dimensional commutative algebra $A$ over a perfect field $k$. When $\text{char}(k) = 0$, the conjecture says this: if $A$ is a subalgebra of a principal ideal algebra $B$, and $\Omega_{A/k}$ injects into $\Omega_{B/k}$, then $A$ is a principal ideal algebra. Here a *principal ideal algebra* is a finite dimensonal commutative $k$-algebra so that every ideal is principal, i.e., of the form $(x)$ for some $x$. To state the conjecture when $\text{char}(k) \neq 0$, we replace 'principal ideal algebra' by 'tame principal ideal algebra'; we will of course define 'tame' and restate the conjecture below.

As the name suggests, this is an Artinian version of a conjecture formulated over 30 years ago by R. Berger in [B]. Berger's Conjecture concerns the coordinate ring $R$ of a reduced curve over a perfect field $k$, and says that $\Omega_{R/k}$ is torsion-free if and only if $R$ is regular. (One direction is classical: if $R$ is regular then $\Omega_{R/k}$ is torsion-free, because it is a projective $R$-module [W, 9.3.14].)

**Main Theorem 0.1 ("ABC ⇒ BC").** *If char(k) = 0, the Artinian Berger Conjecture implies Berger's Conjecture.*

AMS-MOS 1991 Subject Classification: Primary 13N05; Secondary 13E10, 14H20.  
[1] Partially supported by National Security Agency grant MDA 904-92-H-3019  
[2] Partially supported by National Science Foundation grant DMS-9202047  
[3] Our ABC is not to be confused with the *abc conjecture* ("abcc") of number theory.

To see the connection between the two conjectures, let $R$ be the coordinate ring of a singular curve over $k$ with integral closure $S$ and total ring of fractions $F$. Since $\Omega_{S/k}$ is torsion-free, it injects into $\Omega_{F/k} = F \otimes \Omega_{R/k}$. Therefore the torsion submodule of $\Omega_{R/k}$ is the kernel of $\Omega_{R/k} \to \Omega_{S/k}$. It is possible to find a nonzero ideal $I$ of $S$ contained in $R$ so that $B = S/I$ is a principal ideal algebra but $A = R/I$ is not. The torsion submodule of $\Omega_{R/k}$ maps to the kernel $\tau$ of $\Omega_{A/k} \to \Omega_{B/k}$, where we can hope to detect it.

At present, Berger's conjecture is known to be true if $R$ is a complete intersection ([B]), is graded ([S]), has analytically smoothable curve singularities ([Ba]), has multiplicity $e \leq 9$ ([U], [Gu], [I]), or has deviation less than or equal to 3 ([U], [HW]). We refer the reader to the Herzog's survey paper [H] for more details. In this paper, we will use the ABC as a tool to prove that Berger's Conjecture holds in four cases:

1.4 Berger's conjecture holds for seminormal curves (in all characteristics). This folklore result can also be proven using analytic deformation to the graded case, but we know of no literature reference for this result.

1.6 The graded case in characteristic 0. This result is due to Scheja ([S]).

2.13 1-dimensional local rings $(R, M)$ so that $M^3 S \subset R$, where $S$ is the normalization of $R$ and $\mathrm{char}(k) = 0$. This is a new case of Berger's Conjecture.

2.14 Unibranch singularities $(R, M)$ with multiplicity $e < \binom{m}{2}$, where $m = \dim(M/M^2)$ and $\mathrm{char}(k) = 0$. This result is due to Güttes ([Gu], [I]).

Our notion of 'tame' (for tamely ramified) is designed to avoid a pathology in characteristic $p$, namely: "wildly ramified" extensions of $k$ such as $k[s]/(s^p)$ can contain subrings $A$ so that $\Omega_{A/k}$ injects into $\Omega_{B/k}$. Here are two examples which illustrate this pathology.

**Wild Examples 0.2:** a) $A = k[x,y]/(x^3, xy, y^2)$ is isomorphic to the subring $k[s^2, s^3]$ of $B = k[s]/(s^5)$, with $x = s^2$ and $y = s^3$. If $\mathrm{char}(k) = 5$, then $\Omega_{A/k}$ injects into $\Omega_{B/k}$. In fact, a direct computation shows that $\{dx, dy, xdx, ydx\}$ is a $k$-basis of $\Omega_A$, and that it maps to the set $2sds, 3s^2ds, 2s^3ds, 2s^4ds$ of linearly independent elements in $\Omega_B \cong Bds$.

b) Embed $A = k[x,y]/(x^2, xy, y^2)$ in $B = \prod k[s_i]/(s_i^2)$ by setting $x = s_1 + s_3$ and $y = s_2 + s_3$. If $\mathrm{char}(k) = 2$ then $\Omega_{A/k}$ injects into $\Omega_{B/k}$. In this case the vector space $\Omega_{A/k}$ is 5-dimensional on basis $\{dx, dy, xdx, ydy, xdy = -ydx\}$. By inspection, these map to linearly independent elements of the 6-dimensional vector space $\Omega_{B/k} \cong B$.

We now turn to the definition of 'tame' algebra. Since $k$ is a perfect field, every principal ideal algebra $B$ is a finite product of truncated polynomial rings $B_i = K_i[s]/(s^{n_i})$ over finite field extensions $K_i$ of $k$. This classification follows from Wedderburn's Principal Theorem ([Wdb]).

**Definition:** A truncated polynomial ring $B = K[s]/(s^n)$ is said to be *tame* if $K$ is a finite field extension of $k$ and either $\text{char}(K) = 0$, or else $\text{char}(K) = p$ and $p$ does not divide $n$. We say that a principal ideal algebra $B$ is tame if it is the product of tame truncated polynomial rings.

Although the torsion submodule $\tau(R)$ of $\Omega_{R/k}$ only makes sense if $R$ is reduced, we can formulate an Artinian analogue $\tau(A)$. Our definition is motivated by the observation, made above for curves, that $\tau(R)$ is the kernel of $\Omega_{R/k} \to \Omega_{S/k}$. Consider the family $\mathcal{F}$ of submodules of $\Omega_{A/k}$ which arise as the kernel of a map $f_* : \Omega_{A/k} \to \Omega_{B/k}$ induced by an algebra map $f : A \to B$ in which $B$ is a tame principal ideal algebra. This family is closed under intersection, since the product $B \times B'$ of two principal ideal algebras is again a principal ideal algebra, and $\Omega_{B \times B'} = \Omega_B \times \Omega_{B'}$. Because the residue fields of $A$ are tame, $\Omega_{A/k} \in \mathcal{F}$. Since $\Omega_{A/k}$ is an Artinian module, $\mathcal{F}$ has a unique minimal submodule $\tau$.

**Definition:** Let $\tau(A)$ denote the unique minimal submodule of $\mathcal{F}$. Thus $\tau(A)$ is in the kernel of $f_* : \Omega_{A/k} \to \Omega_{B/k}$ for every algebra homomorphism $f : A \to B$ in which $B$ is a tame principal ideal algebra, and equals $\ker(f_*)$ for some $f$. The submodule $\tau(A)$ is natural in $A$; an algebra map $A \to A'$ will induce a map from $\tau(A)$ to $\tau(A')$.

**Artinian Berger Conjecture ("ABC"):** Let $A$ be a finite dimensional commutative algebra over a perfect field $k$. Then:

$$\tau(A) = 0 \Longleftrightarrow A \text{ is a tame principal ideal algebra.}$$

If $A$ is a tame principal ideal algebra, it is clear that $\tau(A) = 0$. If $A$ is a "wild" (not tame) principal ideal algebra, we will see in 2.2 that $\tau(A) \neq 0$. Therefore the ABC is equivalent to the assertion that, if there exists a map to a tame principal ideal algebra $B$ so that $\Omega_A$ injects into $\Omega_B$, then $A$ is a principal ideal algebra. This formulation obviously implies the version of the ABC stated at the outset of this paper; we will show that they are equivalent in 2.4 below.

The paper is organized so as to focus on Berger's conjecture. Section 1 contains the proof of the Main Theorem and a simple proof of Berger's conjecture for seminormal rings. In section 2 we establish the ABC for several classes of local Artinian algebras, including those for which $M^3 = 0$. We reserve section 3 for the proofs of the various technical results about Hochschild homology that we need in sections 1 and 2. In particular, the key result Theorem 1.2 is proven in section three.

*Notation:* All rings in this paper will be commutative algebras over a perfect field $k$. As usual, a *finite* algebra is one that is finite dimensional as a vector space. If $R$ is any algebra, we write $\Omega_R$ for the $R$-module of Kähler differentials $\Omega_{R/k}$ of $R$ over $k$. The terms 'principal ideal algebra' and 'tame principal ideal algebra' were defined above; note that a principal ideal algebra is always a finite algebra.

## §1. Reduction to the Artinian Case.

In this section we prove the main theorem 0.1, modulo some technical results on relative Hochschild homology which we postpone until §3. In order to construct our key commuatative diagram (1.1), we need to review some basic definitions. We refer the reader to [W] for a more detailed discussion of Hochschild homology.

The $R$-module of Kähler differentials $\Omega_{R/k}$ of a commutative algebra $R$ is defined by the following presentation: there is one generator $dx$ for every $x \in R$, with $dx = 0$ if $x \in k$, subject to the usual calculus relations for $d(x + y)$ and $d(xy)$. On the other hand, the Hochschild homology $HH_*(R)$ of $R$ is the homology of the standard Hochschild complex of $R$-modules

$$C_*(R): \quad \cdots R \otimes_k R \otimes_k R \xrightarrow{b} R \otimes_k R \xrightarrow{0} R \to 0.$$
$$b(x \otimes y \otimes z) = xy \otimes z - x \otimes yz + zx \otimes y$$

(See [W, 9.1.1] for a description of the other differentials.) We shall make use of the well-known isomorphism $\Omega_R \cong HH_1(R)$ ([W, 9.2.2]) to fit $\Omega_R$ into a relative homological framework.

If $I$ is an ideal of $R$, we will write $HH_*(R, I)$ for the homology of the kernel $C_*(R, I)$ of the surjection $C_*(R) \to C_*(R/I)$. Our indexing is so that we have a long exact sequence

$$\cdots HH_{n+1}(R/I) \to HH_n(R, I) \to HH_n(R) \to HH_n(R/I) \cdots.$$

One can check directly that $HH_0(R,I) = I$ and that $HH_1(R,I)$ is a quotient of $I \otimes_k R$.

Suppose that $f: R \to S$ is an algebra homomorphism, mapping an ideal $I$ of $R$ isomorphically onto an ideal of $S$ (which by abuse of notation we also write as $I$). We define the double relative groups $HH_*(R, S, I)$ to be the homology of $\mathrm{cone}(f_*)[1]$, the translate of the mapping cone complex of $f_* : C_*(R, I) \to C_*(S, I)$. As in [W, 1.5.2], these fit into a long exact sequence

$$\cdots \xrightarrow{f_*} HH_{n+1}(S,I) \to HH_n(R,S,I) \to HH_n(R,I) \xrightarrow{f_*} HH_n(S,I) \cdots.$$

We will see in 3.6 that $HH_0(R, S, I)$ is isomorphic to $I \otimes_S \Omega_{S/R}$. Assembling all the above data, we form the following commutative diagram with exact rows and an exact column.

(1.1)
$$\begin{array}{ccccccc}
 & HH_1(R,I) & \to & \Omega_R & \to & \Omega_{R/I} & \to 0 \\
 & \downarrow & & \downarrow & & \downarrow & \\
HH_2(S/I) \to & HH_1(S,I) & \to & \Omega_S & \to & \Omega_{S/I} & \to 0 \\
\searrow & \downarrow & & & & & \\
 & I \otimes_S \Omega_{S/R} & & & & &
\end{array}$$

*Diagram Chase 1.1.1:* If the map $HH_2(S/I) \to I \otimes_S \Omega_{S/R}$ is onto (e.g., if $\Omega_{S/R} = 0$), a diagram chase shows that there is an exact "Mayer-Vietoris" sequence

$$\Omega_R \to \Omega_S \oplus \Omega_{R/I} \to \Omega_{S/I} \to 0.$$

Here is the statement of our key technical result. In order to present the flow of ideas more clearly, we shall postpone giving a self-contained proof of this result until §3.

**Theorem 1.2:** *Suppose that an algebra map $f: R \to S$ maps an ideal $I$ of $R$ isomorphically onto an ideal of $S$, that $S$ is locally a principal ideal ring, and that $S/I$ is a finite algebra. If $\mathrm{char}(k) = 0$, then the composite map $HH_2(S/I) \to I \otimes_S \Omega_{S/R}$ in diagram (1.1) is a surjection.*

**Remark:** The cyclic homology and $K$-theory versions of this theorem were proven for number fields in [Wa, A.3], and our proof follows its outline. All these versions are inspired by the Geller-Roberts excision theorem ([GR, 3.1]), which says in effect that $K_2(S, I)$ maps onto $K_1(R, S, I) \cong I \otimes_S \Omega_{S/R}$.

**Corollary 1.3:** *If $\mathrm{char}(k) = 0$, then there is an exact sequence*

$$\Omega_R \to \Omega_S \oplus \Omega_{R/I} \to \Omega_{S/I} \to 0.$$

*Hence if $\Omega_R$ injects into $\Omega_S$ then $\Omega_{R/I}$ also injects into $\Omega_{S/I}$.*

**Proof of the Main Theorem 0.1:** Let $R$ be the coordinate ring of a singular curve over a field $k$ of characteristic 0, and let $S \ne R$ be the normalization of $R$ (the integral closure of $R$ in its total ring of fractions $F$). It is well-known that the conductor ideal $I_0 = \operatorname{ann}_R(S/R)$ has height 1, and that the intersection $J$ of all the singular primes of $R$ is the radical ideal of $I_0$. We set $I = (I_0)^2$, and consider the finite Artinian algebras $A = R/I$ and $B = S/I$. By Nakayama's Lemma. $A = R/I$ is not a principal ideal algebra. On the other hand $B = S/I$ is a principal ideal algebra, since $I$ is a height 1 ideal in a 1-dimensional regular ring. By the ABC, $\ker(\Omega_A \to \Omega_B) \ne 0$. By 1.3, $\ker(\Omega_R \to \Omega_S) \ne 0$. But this is the torsion submodule of $\Omega_{R/k}$. ∎

Here are two more applications of the diagram chase 1.1.1. The first shows that Berger's conjecture holds for seminormal curves in *any* characteristic. The second recovers Scheja's result—Berger's conjecture holds for graded $R$ in characteristic 0.

**Proposition 1.4:** (Folklore) *If $R$ is the coordinate ring of a seminormal singular curve over a perfect field $k$, then $\Omega_R$ has a non-zero torsion submodule.*

**Proof:** Let $M$ be a singular prime ideal of $R$. Since the localization of $\Omega_{R/k}$ at $M$ is $\Omega_{R_M/k}$, we may localize at $M$ to assume that $R$ is local. Since $R$ is seminormal and 1-dimensional, it is well-known that the integral closure $S$ of $R$ satisfies $MS \subseteq R$, and that $K = S/M$ is a finite product of fields (see [T, 1.3]). Therefore $\Omega_{S/R} = \Omega_{K/k} = 0$. By the diagram chase 1.1.1, the torsion submodule of $\Omega_{R/k}$ maps onto the kernel of $(\Omega_{R/I} \to \Omega_{S/I})$ for every ideal $I$, including $I = M^2$. Thus we only need show that $\tau(R/M^2) \ne 0$. Suppose that $\{x_1, \ldots, x_m\}$ is a basis of $M/M^2$. Then $m \ge 2$ and the $\binom{m}{2}$ differentials $x_i dx_j$ with $i < j$ are linearly independent in $\Omega_{R/M^2}$ but vanish in $\Omega_{S/M^2}$. (In fact, we will see in 2.6 below that they form a basis of $\tau(R/M^2)$ when $\operatorname{char}(k) \ne 2$.) ∎

**Proposition 1.5:** *Let $A = k \oplus A_1 \oplus \ldots$ be a graded subalgebra of $B = \prod k[s_i]/(s_i^{n_i})$, $\operatorname{char}(k) = 0$. If $A$ is not a principal ideal algebra then $\Omega_A$ does not inject into $\Omega_B$.*

**Proof:** If the maximal ideal $M$ of $A$ is not principal we can pick homogeneous elements $x \in M_e$, $y \in M_f$ in $M$ which are linearly independent modulo $M^2$. The Euler differential $\omega = ex\, dy - fy\, dx$ in $\Omega_A$ is nonzero, because its image in $\Omega_{A/M^2}$ is $(e+f)x\, dy \ne 0$. But by direct computation $\omega$ vanishes in each factor $\Omega_{k[s]/(s^n)}$ of $\Omega_B$. ∎

**Corollary 1.6:** (Scheja [S]) *Let $R = k \oplus R_1 \oplus \ldots$ be a graded reduced 1-dimensional algebra of finite type over a field $k$ of characteristic zero. If $R \ne k[t]$ then $(\Omega_R)_{tor} \ne 0$.*

**Proof:** In this case the normalization $S$ is also graded. By base change we can assume all the residue fields of $S$ are $k$. If the homogeneous maximal ideal $M$ of $R$ is principal, then $R = k[t]$. Let $I$ be the ideal $J^2$, $J = \text{ann}_R(S/R)$; as $M$ is the only associated prime of the graded $S/R$, $I$ has height 1. We now quote Proposition 1.5 with $A = R/I$ and $B = S/I$ and use 1.3 to complete the proof. ∎

## §2. Evidence for the truth of ABC.

In this section we develop some tools for detecting $\tau(A)$, and show that the ABC holds for several classes of finite algebras $A$. We also show that $\tau(A)$ is a subspace of the cyclic homology group $HC_1(A) = \Omega_A/dA$ when $\text{char}(k) = 0$ and describe the quotient space. We start with some reductions. The next lemma says that we can always assume, without loss of generality, that $k$ is algebraically closed and that $A$ is local.

**Lemma 2.0:** *Let $A$ be a finite commutative $k$-algebra.*
  *(a) If $K$ is a finite field extension of $k$, then under the canonical isomorphism $\Omega_{A \otimes_k K} \cong \Omega_A \otimes_k K$ we have $\tau(A \otimes_k K) \cong \tau(A) \otimes_k K$.*
  *(b) If $A = A_1 \times \cdots \times A_n$ then the decomposition $\Omega_A = \oplus \Omega_{A_i}$ induces a decomposition $\tau(A) = \oplus \tau(A_i)$.*

**Proof:** Let $f: A \to B$ be a homomorphism with $\tau(A) = \ker(f_*)$. Then $A \otimes K \to B \otimes K$ induces a map $\Omega_{A \otimes K} \to \Omega_{B \otimes K}$ with kernel $\tau(A) \otimes K$. Hence $\tau(A \otimes K) \subseteq \tau(A) \otimes K$. But $\tau(A) \subseteq \tau(A \otimes K)$ by naturality, so we have equality in (a). If $A = \prod A_i$ then by naturality each $\tau(A_i)$ (and hence the subspace $\oplus \tau(A_i)$ of $\Omega_A$) lies in $\tau(A)$. To see equality, choose $f_i: A_i \to B_i$ with $\tau(A_i) = \ker(f_{i*})$; $\oplus \tau(A_i)$ is the kernel of $(\prod f_i)_*: \Omega_A \to \Omega_{B_1 \times \cdots \times B_n}$. ∎

**Lemma 2.1:** *Let $B$ be a tame principal ideal algebra. If an element $x \in B$ satisfies $x^i = 0$, then $x^{i-1}dx = 0$ in $\Omega_{B/k}$.*

**Proof:** This is clear if $i \neq 0$ in $k$: $x^{i-1}dx = d(x^i/i)$. Suppose that $i = 0$ in $k$ and that $s$ is a parameter of some truncated polynomial ring which is a factor of $B$. If $s^i = 0$ then $s^{i-1} = 0$, and $s^{i-1}ds = \frac{-s}{(i-1)} d(s^{i-1}) = 0$. In general, if the leading term of $x$ is $\alpha s^e$, then $s^{ei} = 0$ and hence $x^{i-1}dx = (\alpha^i e + \ldots)s^{ei-1}ds = 0$. ∎

**Corollary 2.2:** *If $A$ is the wild principal ideal algebra $K[x]/(x^{np})$ then $\tau(A) = K\, x^{np-1}dx$. Thus $\tau(A) \neq 0$ for every wildly ramified principal ideal algebra $A$.*

**Proof:** Embed $A$ in $K[x]/(s^{np(p+1)-1})$ by $x \mapsto s^{p+1}$. Then $x^{i-1}dx$ maps to $s^{i(p+1)-1}ds$, which is nonzero for $i < np$. Thus $\tau(A) \subseteq K\, x^{np-1}dx$. But $x^{np-1}dx \in \tau(A)$ by 2.1. ∎

In order to simplify our computations we next show that we may restrict our attention to subalgebras $A$ of tame principal ideal algebras.

**Definition:** We say that a finite $k$-algebra $A$ is *embeddable* if it is isomorphic to a subalgebra of some tame principal ideal algebra. For example, if $M^2 = 0$ then $A$ is embeddable into a product of $\dim_k(M)$ truncated polynomial rings. Wild principal ideal algebras may also be embeddable, as the proof of 2.2 shows.

**Lemma 2.3:** *Every finite $k$-algebra $A$ has a maximal embeddable quotient $\bar{A}$. Moreover, $\bar{A}$ is not a principal ideal algebra unless $A$ is, and $\tau(A)$ maps onto $\tau(\bar{A})$.*

**Proof:** Let $A$ be a finite $k$-algebra, and consider the family $\mathcal{F}$ of all ideals $I$ of $A$ so that $A/I$ is embeddable. This is not the empty set, because we have seen that $M^2 \in \mathcal{F}$ for every maximal ideal $M$ of $A$. If $I_1$ and $I_2$ are in $\mathcal{F}$, then $I_1 \cap I_2 \in \mathcal{F}$, because $A/(I_1 \cap I_2) \hookrightarrow A/I_1 \times A/I_2 \hookrightarrow B_1 \times B_2$. By the descending chain condition, there is a unique minimal ideal $I_{min}$ in $\mathcal{F}$. By construction, $\bar{A} = A/I_{min}$ is embeddable. Moreover, any map $f$ from $A$ to a tame principal ideal algebra $B$ must factor through $\bar{A}$, because $\ker(f) \in \mathcal{F}$. Thus there is a one-one correspondence between homomorphisms $f : A \to B$ and homomorphisms $\bar{f} : \bar{A} \to B$. It follows that the surjection $\Omega_A \to \Omega_{\bar{A}}$ maps $\tau(A)$ onto $\tau(\bar{A})$. Finally, if $A$ is not a principal ideal algebra, then some maximal ideal $M$ of $A$ is not principal. Since $I_{min} \subseteq M^2$, Nakayama's lemma implies that $A/I_{min}$ is not a principal ideal algebra. ∎

**Criterion 2.3.1:** In a truncated polynomial ring any solution to $x^2 = y^2 = 0$ must satisfy $xy = 0$. Equivalently, if $xy \neq 0$ then either $x^2 \neq 0$ or $y^2 \neq 0$. This gives a simple test for non-embeddability. For example, the algebra $A = k[x,y]/(x^2, y^2)$ is not embeddable. In this case $M^2 = (xy)$, so the maximal embeddable quotient is $\bar{A} = A/M^2 = k[x,y]/(x^2, xy, y^2)$.

**Proposition 2.4:** *The Artinian Berger Conjecture is equivalent to the assertion:*

*if $A$ is a subalgebra of a tame principal ideal algebra $B$ so that $\Omega_A \hookrightarrow \Omega_B$,*

*then $A$ is a tame principal ideal algebra.*

**Proof:** Since $\Omega_A \hookrightarrow \Omega_B$ implies that $\tau(A) = 0$, the ABC implies the displayed assertion. If $A$ is a principal ideal algebra, then 2.2 shows that $\tau(A) = 0$ iff $A$ is tame. We must show that if $A$ is not a principal ideal algebra then $\tau(A) \neq 0$. By Lemma 2.3, we may replace $A$ by $\bar{A}$ to suppose that there exists an embedding of $A$ into some tame principal ideal algebra $B'$. By construction, $\tau(A)$ is the kernel of $(\Omega_A \to \Omega_{B''})$ for some homomorphism $A \to B''$ with $B''$ tame. Then $A$ embeds into $B = B' \times B''$, and the minimality of $\tau(A)$ implies that $\tau(A) = \ker(\Omega_A \to \Omega_B)$. Since $A$ is not a principal ideal algebra, the displayed assertion implies that $\tau(A) \neq 0$, as desired. ∎

Here is a criterion for an element of $\Omega_A$ to lie in $\tau(A)$; Lemma 2.1 is a special case.

**Lemma 2.5:** *Let $A$ be a finite local $k$-algebra. Suppose that $x, y \in A$ satisfy $xy = 0$. Then $x\, dy \in \tau(A)$.*

**Proof:** It suffices to show that $x\, dy$ vanishes in $\Omega_B$ for every map $g\colon A \to B$ in which $B = k[s]/(s^n)$ is a tame truncated polynomial ring. There are nonzero constants $\alpha, \beta \in k$ so that $g(x) = \alpha s^e + u s^{e+1}$ and $g(y) = \beta s^f + v s^{f+1}$ with $u, v \in B$. Since $xy = 0$, we have $s^{e+f} = 0$ in $B$. Then the image of $xdy$ in $\Omega_B$ is $(\alpha s^e + u s^{e+1}) d(\beta s^f + v s^{f+1}) = \alpha \beta f s^{e+f-1} ds$, which is zero by Lemma 2.1. ∎

**Seminormal Example 2.5.1:** Let $A$ be the subalgebra of $B = \prod k[s_i]/(s_i^{n_i})$ generated by $\{s_1, ..., s_m\}$, i.e., $A \cong k[x_1, \ldots, x_m]/(s_i^{n_i}, s_i s_j \text{ for } i \neq j)$. A straightforward calculation shows that the kernel of $\Omega_A \to \Omega_B$ has for a basis the set of all $s_i\, ds_j$, $i < j$. If $B$ is tame, 2.5 shows that they form a basis for $\tau(A)$.

**Proposition 2.6:** *Suppose that $A = K \oplus M$ with $M^2 = 0$, where $K$ is a finite extension of $k$, and let $\{\bar{x}_1, \ldots, \bar{x}_m\}$ be a basis of $M/M^2$. Then $\Omega_A/\tau(A) \cong M$, and:*

*If $\operatorname{char}(k) \neq 2$ then $\dim \tau(A) = \binom{m}{2}$, and the $x_i\, dx_j$, $i < j$, form a basis of $\tau(A)$.*

*If $\operatorname{char}(k) = 2$ then $\dim \tau(A) = \binom{m+1}{2}$, and the $x_i\, dx_j$, $i \leq j$, form a basis of $\tau(A)$.*

**Proof:** Since $\Omega_{A/k} = \Omega_{A/K}$, we may assume that $K = k$. We use the calculation of $\Omega_A$ given below in 3.2. Suppose first that $\operatorname{char}(k) \neq 2$. Map $A$ to the tame $B = \prod k[x_i]/(x_i^2)$ in the obvious way. The map from $\Omega_A$ to $\Omega_B \cong M$ is the map $\mu$ of 3.2. Hence $\tau(A)$ lies in $\ker(\mu)$, which by 3.2 is isomorphic to $\Lambda^2 M$ and has the $x_i\, dx_j$ as a basis. But $x_i dx_j \in \tau(A)$ by 2.5.

If $\operatorname{char}(k) = 2$, we map $A$ to the tame $B = \prod k[s_i]/(s_i^3)$ by sending $x_i$ to $s_i^2$. Again, the image of $\Omega_A$ in $\Omega_B$ is isomorphic to $M$, and by 3.2 the kernel is isomorphic to the vector space $\tilde{\Lambda}^2 M$ spanned by the $x_i\, dx_j$, $i \leq j$. These differentials are in $\tau(A)$ by 2.5. ∎

**Corollary 2.7:** Let $A$ be a finite local $k$-algebra with maximal ideal $M$. Suppose $xy = 0$ for two elements $x, y \in M$ which are linearly independent mod $M^2$. Then $x\,dy$ is a nonzero element of $\tau(A)$.

**Proof:** We have $x\,dy \in \tau(A)$ by 2.5, and it is nonzero because $0 \neq \bar{x}\,d\bar{y} \in \Omega_{A/M^2}$. ∎

**Application 2.8:** The *socle* $I = \operatorname{ann}_A(M)$ is a nonzero ideal in any Artinian local ring $A$ (except a field), and there is a map $I \otimes_A M \to \Omega_{A/k}$ sending $x \otimes y$ to $xdy$. The image $V_A$ of this map is a submodule of $\tau(A)$ by 2.5, so $A$ will satisfy the ABC whenever $V_A \neq 0$. This happens, for example, in the following two cases.

   a) If the socle is not contained in $M^2$, and $\dim M/M^2 \neq 1$, then $\tau(A) \neq 0$.

   b) If $A = k[x_1, ..., x_m]/(x_1, ..., x_m)^3$ and $\operatorname{char}(k) \neq 3$, the socle $I = M^2$ has dimension $\binom{m+1}{2}$. Since $\Omega_A$ is the quotient of the free $A$-module on $\{dx_i\}$ by the $\binom{m+2}{3}$ relations $d(x_i x_j x_k) = 0$, it is straightforward to see that $\tau(A) = V_A$, and that its dimension is $2\binom{m+1}{3} = 2\binom{m}{3} + 2\binom{m}{2} = m\binom{m+1}{2} - \binom{m+2}{3}$.

   Our next result establishes the ABC for the family of algebras

$$A_m = k[x, y]/(x^m, x^2y, y^2).$$

**Proposition 2.9:** Let $A$ be a finite local algebra which is not a principal ideal algebra. If there exists a $y \in A$ so that $A/yA$ is a principal ideal algebra and that $\dim yA \leq 2$, then $\tau(A) \neq 0$.

**Proof:** We may assume the residue field of $A$ is $k$. Since $A/yA$ is local, $A/yA \cong k[x]/(x^m)$ for some $x \in A$ and $m \geq 2$. Since $M \neq xA$, Nakayama's Lemma gives $y \notin M^2$. If $xy = 0$ we are done by 2.7, so we may assume that $xy \neq 0$. Thus $\{y, xy\}$ forms a basis of $yA$. Write $y^2 = \alpha xy$ and $x^m = \beta xy$. Since $y(y - \alpha x) = x(\beta y - x^{m-1}) = 0$ in $A$, we are again done by 2.7 unless $\alpha = \beta = 0$, i.e., $A = A_m$ for some $m \geq 2$. If $m = 2$ then $A_2$ is not embeddable; we saw in 2.3.1 that $\bar{A}_2 = k \oplus (M/M^2)$. Since $\tau(\bar{A}_2) \neq 0$ by 2.6, we have $\tau(A) \neq 0$ by 2.3.

We have reduced to the case $A = A_m$, $m \geq 3$. A direct calculation shows that $\{x^i dx | \ i = 0, 1, ..., m-2\} \cup \{dy, xdy, ydx, xydx\}$ is a linearly independent subset of $\Omega_{A_m}$, and forms a basis if $\operatorname{char}(k) \nmid m$. In this case the socle contains $xy$, so $xy\,dx \in \tau(A_m)$ by 2.8. Hence $\tau(A_m)$ is nonzero. ∎

**Remark 2.9.1:** The kernel of $\Omega_{A_3} \to \Omega_B$ depends upon the choice of $B$, even when we restrict to graded algebra maps. If we embed $A_3$ into $B_1 = k[s]/(s^6)$ by setting $x = s^2$, $y = s^3$, then $w_1 = 2xdy - 3ydx$ and $xydx$ form a basis of $\ker(\Omega_{A_3} \to \Omega_B)$. But if we embed $A_3$ into $B_2 = k[t]/(t^9)$ by setting $x = t^3$, $y = t^5$, then $w_2 = 3xdy - 5ydx$ and $xydx$ form a basis of $\ker(\Omega_{A_3} \to \Omega_B)$. This shows that $\tau(A)$ is generated by $xydx$.

Our next result handles subrings $A$ of a product $B = \prod k[s_i]/(s_i^3)$. The following definition will be useful.

**Definition 2.10.0:** Given an algebra map $\pi: A \to B = k[s]/(s^n)$, we define the function $\nu: A \to \{0, 1, ..., n-1, \infty\}$ by: $\nu(a) = e$ when $\pi(a) = \alpha s^e + ws^{e+1} \in B_i$ for some $0 \neq \alpha \in k$ and $w \in B$; we set $\nu(a) = \infty$ if $\pi(a) = 0$.

**Proposition 2.10:** Let $A$ be a finite local $k$-algebra which is not a principal ideal algebra. Suppose that $\dim \pi(M^2) \leq 1$ for every map $\pi : A \to k[s]/(s^n)$ into a tame truncated polynomial algebra. Then $\tau(A) \neq 0$.

**Proof:** We first show that $\omega = x\,dy - y\,dx \in \Omega_A$ belongs to $\tau(A)$ whenever $x, y \in M$. For this it suffices to show that $\pi(\omega) = 0$ in $\Omega_B$ for every tame $B = k[s]/(s^n)$ and every map $\pi: A \to B$. Proceeding as in the proof of 2.5, we write $\pi(x) = \alpha s^e + us^{e+1}$ and $\pi(y) = \beta s^f + vs^{f+1}$, where $0 \neq \alpha, \beta \in k$ and $u, v \in B$. If $\pi(xy) = 0$ then (as in 2.5) we have $\pi(x\,dy) = \pi(y\,dx) = 0$, and hence $\pi(\omega) = 0$ in $\Omega_B$. If $\pi(xy) \neq 0$, our hypothesis forces $\pi(x^2)$ and $\pi(y^2)$ to be scalar multiples of $\pi(xy)$, and hence $e = f$. Write $\pi(x^2) = \lambda \pi(xy)$; we see that

$$s^{2e}(\alpha + us)^2 = \lambda s^{2e}(\alpha + us)(\beta + vs).$$

Since $\alpha + us$ is invertible in $B$ this yields $s^e \pi(x) = \lambda s^e \pi(y)$. Thus the polynomial $h = \pi(\lambda y - x)$ satisfies $s^e h = 0$. By Lemma 2.5, $\pi(y)\,dh = h\,d\pi(y) = 0$, whence in $\Omega_B$ we have $\pi(x\,dy) = \lambda \pi(y\,dy) = \pi(y\,dx)$, or $\pi(\omega) = 0$. Therefore $\omega \in \tau(A)$.

If $\text{char}(k) \neq 2$ this proves that $\tau(A) \neq 0$, because if $x, y$ are linearly independent mod $M^2$ then $\omega \neq 0$. Indeed, the image of $\omega$ in $\Omega_{A/M^2}$ is $2\bar{x}\,d\bar{y}$, which is nonzero by 2.6.

It remains to prove the result when $\text{char}(k) = 2$. By Lemma 2.3, we may replace $A$ by $\bar{A}$ if necessary in order to assume that $A$ is embeddable into some tame principal ideal algebra $B = \prod B_i$, $B_i = k[s_i]/(s_i^{n_i})$. By 2.7 we may assume that $xy \neq 0$ for every $x, y \in M$ which are linearly independent mod $M^2$. By Lemma 2.10.1 below, we can find $x_1, ..., x_m$ mapping to a basis of $M/M^2$ so that $A$ has the presentation (2.10.2). Since $0 \neq x_i x_j \in V$ for every $i < j$, the final assertion in 2.10.1 is that each $\omega = d(x_i x_j) \neq 0$ in $\Omega_A$. Since we have seen above that $\omega \in \tau(A)$, we are done. ∎

**Lemma 2.10.1:** *Suppose that $A$ is an embeddable algebra satisfying the hypotheses of 2.10. Suppose moreover that $xy \neq 0$ for every $x, y \in M$ which are linearly independent mod $M^2$. Then there exist $x_i \in M$ and $c_{ijk} \in k$ so that*

$$(2.10.2) \quad A \cong k[x_1, ..., x_m]/I, \quad I = (x_1, ..., x_m)^3 + \left(\sum_{i<j} c_{ijk} x_i x_j, \quad k = 1, ..., N\right).$$

*If $V$ denotes the subspace of $M^2$ spanned by $\{x_i x_j \mid i < j\}$ then $d: V \to \Omega_A$ is an injection.*

**Proof:** We first claim that $x^2 \neq 0$ for every $x \in M - M^2$. To see this, choose $y$ with $xy \neq 0$ and find a projection $f: A \to k[s]/(s^n)$ with $f(xy) \neq 0$. Then $\nu(xy) = \nu(x) + \nu(y) < \infty$. If $\nu(x) < \nu(y)$ then $\nu(x^2) = 2\nu(x) < \nu(xy)$, contradicting the assumption that $\dim f(M^2) = 1$. Similarly, we cannot have $\nu(y) < \nu(x)$. Hence $\nu(y) = \nu(x)$ and $\nu(x^2) = \nu(xy) < \infty$, which implies that $x^2 \neq 0$.

Second, we shall see that $M^3 = 0$. Since $A$ is embeddable, it suffices to show that $f(M^3) = 0$ for every map $f: A \to k[s]/(s^n)$. Choose $x \in M$ of minimum valuation $e > 0$. Suppose that $f(M^3) \neq 0$. Since $f(M^3) \subseteq s^{3e} S$ we have $\nu(x^2) = 2e < 3e = \nu(x^3) < \infty$. But this contradicts the assumption that $\dim f(M^2) \leq 1$.

Next we observe that for each $x \in M - M^2$ there is a surjection $\pi: A \to k[s]/(s^3)$ sending $x$ to $s$. To see this, choose a projection $f: A \to k[t]/(t^n)$ in which $f(x^2) \neq 0$. Since $\dim f(M^2) = 1$, $e = \nu(x) = \min\{\nu(m) \mid m \in M\}$ and we can choose $t_2, ..., t_m \in M$ with $\nu(t_i) > e$ for all $i$ so that $\{x, t_2, ..., t_m\}$ maps to a basis of $M/M^2$. Since $\nu(t_i) > e$, we have $f(t_i M) = 0$ for all $i$. Hence for $I = \ker(f) + (t_2, ..., t_m)$ we have $A/I \cong k[x]/(x^3)$, as claimed.

Proceeding inductively for $m \leq \dim(M/M^2)$, we construct a sequence $x_1, ..., x_m$ of elements in $M$ which are linearly independent mod $M^2$ and a map $f_m: A \to \prod_{i=1}^m k[s_i]/(s_i^3)$ with $f_m(x_i) = s_i$ for all $i$. For the inductive step, choose $x_m \in \ker(f_{m-1})$ and construct $\pi: A \to k[s_m]/(s_m^3)$ as above; if $\pi(x_i) = \alpha_i s_m + \beta_i s_m^2$ we replace $x_i$ by $x_i - \alpha_i x_m - \beta_i x_m^2$ to get $\pi(x_i) = 0$ and arrange that $f_m = f_{m-1} \times \pi$ satisfies $f_m(x_i) = s_i$.

When $m = \dim(M/M^2)$ the sequence $x_1, ..., x_m$ maps to a basis of $M/M^2$, so $A$ is a quotient of $k[x_1, ..., x_m]/(x_1, ..., x_m)^3$. If any quadratic relation $\sum c_{ij} x_i x_j = 0$ holds in $A$ then by applying $f_m$ we see that $c_{ii} = 0$ for all $i$. This gives the presentation of $A$.

Now $\Omega_A$ is the quotient of the free $A$-module on the $dx_i$ by relations $d(xyz) = 0$ of degree three and the quadratic relations $\sum c_{ijk} d(x_i x_j) = 0$. Hence $d(V)$ is the vector space generated by the symbols $d(x_i x_j)$ with the quadratic relations. It follows that the vector spaces $V$ and $d(V)$ are isomorphic. ∎

We now prove that the ABC holds whenever $M^3 = 0$. (The case $M^2 = 0$ is covered by Proposition 2.6.) This will give us a new case of Berger's conjecture.

**Theorem 2.11:** *Let $(A, M)$ be a finite local $k$-algebra satisfying $M^3 = 0$. If $A$ is not a principal ideal algebra then $\tau(A) \neq 0$.*

**Proof:** By theorem 2.10 we may assume that there exists a map $\pi \colon A \to k[s]/(s^n)$ with $\dim \pi(M^2) \geq 2$. Let $0 < e < f$ be the lowest values in the set $\nu(A)$, and choose $x, y \in M$ with $\nu(x) = e$, $\nu(y) = f$, and $y \notin M^2$. Then choose $z_3, ..., z_m \in M$ with $\nu(z_j) > f$ so that $\{x, y, z_3, ..., z_m\}$ maps to a basis of $M/M^2$. Setting $I = \{a \in M^2 \mid \nu(a) > e + f\}$, we have arranged that $M z_j \subseteq I$ for $j = 3, ..., m$, and $\dim(M^2/I) = 2$. Thus $A/(I + (z_3, ..., z_m)A)$ is isomorphic to the ring $A_3 = k[x, y]/(x^3, y^2, x^2 y)$ of 2.9. With these choices, $xy\, dx \in \tau(A)$ by 2.5, and it is nonzero because it maps to $xy\, dx \neq 0$ in $\Omega_{A_3}$. ∎

**Remark 2.11.1:** In summary, if $(A, M)$ is embeddable and $M^3 = 0$ then either:
  1) There exists a map $\pi \colon A \to k[s]/(s^n)$ with $\dim \pi(M^2) \geq 2$, or
  2) The complement of 1).

In case 1), there exist $x, y \in M$ such that $0 \neq xy\, dx \in \tau(A)$; in case 2) $\omega = x\, dy - y\, dx \in \tau(A)$ for all $x, y \in M$, and $\omega \neq 0$ whenever $x, y$ are linearly independent mod $M^2$.

If we use Theorem 2.11 in conjunction with Corollary 1.3 and the Main Theorem 0.1, we get the following two corollaries.

**Corollary 2.12:** *Suppose that $\mathrm{char}(k) = 0$. Let $A$ be a finite local $k$-algebra with maximal ideal $M$. Suppose that $A$ embeds in a principal ideal algebra $B$ so that $M^3 B$ lies in $A$. Then $\Omega_A \to \Omega_B$ is not injective.*

**Theorem 2.13 (Berger's Conjecture if $M^3$ is in the conductor):** *Let $\mathrm{char}(k) = 0$. Suppose that $(R, M)$ is the local ring of a curve, and that $S$ is the normalization of $R$. If $M^3$ is contained in the conductor of $S/R$ (i.e., if $M^3 S \subset R$), then Berger's conjecture holds for $R$.*

As another application of Lemma 2.5, we give a new proof of a theorem of Güttes ([Gu],[I]) relating the embedding dimension $m = \dim(M/M^2)$ to the multiplicity $e$ of a unibranch singularity.

**Proposition 2.14:** *Let $A$ be a subalgebra of $B = k[s]/(s^n)$. Assume that the maximal ideal $M$ of $A$ has embedding dimension $m = \dim(M/M^2)$, and let $e \geq 1$ be maximal so that $M \subseteq s^e B$. If $e < \binom{m}{2}$ then the kernel of $\Omega_A \to \Omega_B$ is nonzero.*

**Proof:** Choose $x \in M$ not in $s^{e+1}B$, and set $I = M^2 + xA$. Since $\Omega_A/x\Omega_A$ surjects onto $\Omega_{A/I}$, the calculation of 3.2 (cited in 2.6) shows that, as a vector space,

$$\dim(\Omega_A/x\Omega_A) \geq \dim(\Omega_{A/I}) \geq \binom{m-1}{2} + (m-1) = \binom{m}{2}.$$

On the other hand, the $k[x]$-module $\Omega_B$ has $e$ generators. Hence, by the general theory of finitely generated torsion modules over the principal ideal domain $k[x]$, any $k[x]$-submodule $L$ of $\Omega_B \cong B/(ns^{n-1})$ can have at most $e$ generators, i.e., $e \geq \dim(L/xL)$. If $\Omega_A \to \Omega_B$ were injective, this would yield the inequality

$$e \geq \dim(\Omega_A/x\Omega_A) \geq \binom{m}{2}. \quad \blacksquare$$

**Corollary 2.15:** (Güttes ([Gu, Satz 5])) *Assume $\mathrm{char}(k) = 0$. If $(R, M)$ is the local ring of a unibranch curve whose multiplicity $e$ and embedding dimension $m = \dim(M/M^2)$ satisfy $e < \binom{m}{2}$, then the torsion submodule $\tau(R)$ of $\Omega_R$ is nonzero.*

**Proof:** The unibranch hypothesis means that the integral closure $S$ of $R$ is local, say with parameter $s$. The multiplicity of $R$ is the largest integer $e$ so that $M \subseteq s^e S$, i.e., the integer such that $MS = s^e S$; we have $e(M, R) = e(M, S) = \dim S/MS = e$ by [BH, 4.6.9]. Choose an ideal of the form $I = s^n S$ contained in $M^2$. Then $A = R/I \subset B = S/I$ satisfy the hypotheses of 2.14, and $\tau(R)$ maps onto $\ker(\Omega_A \to \Omega_B)$ by Corollary 1.4. $\blacksquare$

We conclude this section with a final piece of evidence for the ABC. Recall that the derivative $d : A \to \Omega_A$ is a $k$-linear map whose kernel is the de Rham cohomology group $H^0_{dR}(A)$ and whose cokernel $\Omega_A/dA$ is the cyclic homology group $HC_1(A)$. The following result shows that $HC_1(A)$ is an upper bound for $\tau(A)$.

**Proposition 2.16:** Let $A$ be a finite algebra over a field $k$ of characteristic 0. Then

i) $A$ is a principal ideal algebra $\iff$ $HC_1(A) = 0$.

ii) If $A$ is a subalgebra of a principal ideal algebra $B$ then $H^0_{dR}(A) \cong A_{red}$. Furthermore integration ($\int$) with respect to the parameters $s_i$ of $B$ defines an exact sequence of vector spaces:

$$0 \to \ker(\Omega_A \to \Omega_B) \to HC_1(A) \xrightarrow{\int} B/A \to \Omega_B/\Omega_A \to 0.$$

**Proof:** Suppose first that $A$ is a principal ideal algebra. The usual integration formulas applied to its factors $K_i[s]/(s^{n_i})$ show that $H^0_{dR}(A) = \prod K_i \cong A_{red}$, and $HC_1(A) = 0$. If $A$ is not a principal ideal algebra, some maximal ideal $M$ is not principal. But then $HC_1(A)$ maps onto $HC_1(A/M^2)$, and part (i) follows from the consequence $HC_1(A/M^2) \cong \Lambda^2(M/M^2) \neq 0$ of the calculation of $\Omega_{A/M^2}$ in 3.2.

Now suppose that $A$ is a subalgebra of a principal ideal algebra $B$. The following diagram commutes and has exact rows by the definition of $H^0_{dR}$ and $HC_1$.

$$\begin{array}{ccccccccccc}
0 & \to & H^0_{dR}(A) & \to & A & \xrightarrow{d} & \Omega_A & \to & HC_1(A) & \to & 0 \\
& & \downarrow & & \downarrow & & \downarrow & & \downarrow & & \\
0 & \to & H^0_{dR}(B) & \to & B & \xrightarrow{d} & \Omega_B & \to & HC_1(B) & \to & 0
\end{array}$$

Since $HC_1(B) = 0$, integration gives a map $\Omega_B \to B$ and hence a map from $HC_1(A)$ to $B/A$. The nilradical nil$(A)$ of $A$ lies inside the nilradical of $B$, which injects into $\Omega_B$ by part (i). Thus nil$(A)$ injects into $\Omega_A$. Since $A$ is a product of Artin local rings, each containing a coefficient field $K_i$, we have $\prod K_i \cong A_{red}$ and a vector space decomposition $A = \prod K_i \oplus$ nil$(A)$. Since each $K_i$ is separable over $k$, we have $dK_i = 0$ for all $i$. Hence $\prod K_i = H^0_{dR}(A)$. The exact sequence follows from the snake lemma. ∎

## §3. Relative Hochschild Homology

This is the technical section in which we prove Theorem 1.2 as well as several other assertions about the Hochschild homology and the relative homologies, $HH_1(R, I)$ and $HH_0(R, S, I)$, used in this paper. All rings in this section will be commutative algebras over a fixed field $k$, which need not be perfect, and we shall adopt the notation that $\otimes$ denotes $\otimes_k$.

We start with a description of the relative term $HH_1(R,I)$. Since the kernel of $R \otimes R \to R/I \otimes R/I$ is $R \otimes I + I \otimes R$, it follows from the definition (given in §1) that $HH_1(R,I)$ is the cokernel of the Hochschild boundary map

(3.0) $$b: R \otimes R \otimes I + R \otimes I \otimes R + I \otimes R \otimes R \longrightarrow R \otimes I + I \otimes R$$
$$b(x \otimes y \otimes z) = xy \otimes z - x \otimes yz + zx \otimes y$$

The submodule $I \otimes k$ of $I \otimes R$ is the image under $b$ of the 'degenerate' terms (those with $y$ or $z$ in $k$), so it maps to zero in $HH_1(R,I)$. Thus if $M$ is a maximal ideal of $R$ with $R/M = k$ we may ignore degeneracies by replacing the source and target of $b$ by $R \otimes M \otimes I + R \otimes I \otimes M + I \otimes M \otimes M$ and $R \otimes I + I \otimes M$, respectively.

**Proposition 3.1:** *Suppose that $M$ is a maximal ideal of $R$ with $R/M = k$. If $I \subseteq M$ is an ideal with $IM = 0$, then there is an exact sequence of $R$-modules*

$$I \otimes_k I \xrightarrow{\eta} I \otimes_k (M/M^2) \xrightarrow{\iota} HH_1(R,I) \xrightarrow{\mu} I \to 0$$

*where $\eta(x \otimes y) = x \otimes y + y \otimes x$, $\iota(x \otimes \bar{y}) = x \otimes y$, and $\mu(x \otimes y) = xy$.*

**Proof:** The $R$-linear surjection $\mu : R \otimes I + I \otimes M \to I$ defined by $\mu(x \otimes y) = xy$ has kernel $M \otimes I + I \otimes M$. Since $\mu b$ vanishes on $R \otimes M \otimes I + R \otimes I \otimes M + I \otimes M \otimes M$, $\mu$ induces a well-defined map from $HH_1(R,I)$ onto $I$. If $x \in I$ and $z \in M$ then $b(x \otimes y \otimes z) = xy \otimes z - x \otimes yz$, so $I \otimes M \to HH_1(R,I)$ factors through a map $\iota$ from $I \otimes_R M = I \otimes_R (M/M^2) = I \otimes_k M/M^2$ to $HH_1(R,I)$. This finishes the construction of the sequence. The following formula shows that the sequence is exact at $HH_1(R,I)$, and that $\iota\eta = 0$.

$$b(1 \otimes x \otimes y) = x \otimes y + y \otimes x, \quad x \in I, \ y \in M.$$

To establish exactness at $I \otimes (M/M^2)$, observe that the formula

$$f(x \otimes y) = -f(y \otimes x) = x \otimes y, \quad x \in I, \ y \in M$$

determines a $k$-linear map $f$ from $M \otimes I + I \otimes M$ to $\operatorname{coker}(\eta)$. If $x \in I$ and $y, z \in M$, then

$$fb(x \otimes y \otimes z) = fb(y \otimes x \otimes z) = fb(y \otimes z \otimes x) = -x \otimes yz = 0.$$

Thus $f$ induces a well-defined map from $\ker(HH_1(R,I) \to I)$ to $\operatorname{coker}(\eta)$. By definition, the composition $f\iota$ is the natural projection from $I \otimes (M/M^2)$ to $\operatorname{coker}(\eta)$. This means that the sequence is exact at $I \otimes (M/M^2)$. ∎

**Remark 3.1.1:** If we relax the hypothesis about the residue field to say that $K = R/M$ is a finite separable extension of $k$, the exact sequence becomes

$$I \otimes_K I \xrightarrow{\eta} I \otimes_K (M/M^2) \xrightarrow{\iota} HH_1(R, I) \xrightarrow{\mu} I \to 0.$$

This follows from étale descent ([WG]) applied to a version of 3.1 for localizations of $R \otimes K$ and $I \otimes K$ over the field $K$. Since we will not need this result, we omit the details.

When $R/M = k$, there is an isomorphism $\Omega_R \cong HH_1(R) \cong HH_1(R, M)$. This yields the following corollary, in which $\tilde{\Lambda}^2 M$ denotes $M \otimes_k M/(\{x \otimes y + y \otimes x \mid x, y \in M\})$. Note that if $\text{char}(k) \neq 2$ then $\tilde{\Lambda}^2 M$ is the usual exterior product $\Lambda^2 M$.

**Corollary 3.2:** If $M^2 = 0$ and $R/M = k$, there is a short exact sequence of $R$-modules

$$0 \to \tilde{\Lambda}^2 M \xrightarrow{\iota} \Omega_R \xrightarrow{\mu} M \to 0.$$

**Corollary 3.3:** If $IM = 0$, $R/M = k$, and $I \subseteq M^2$, then there is a short exact sequence

$$0 \to I \otimes_k (M/M^2) \xrightarrow{\iota} HH_1(R, I) \xrightarrow{\mu} I \to 0.$$

**Corollary 3.4:** Suppose that $R/M = k$, $I \cong kt$, and $IM = 0$. Then $HH_1(R, I)$ is isomorphic to the $R$-module $R/(M^2, 2t)$ on generator $dt = 1 \otimes t$.

Next we turn to the double-relative group $HH_0(R, S, I)$. We will work in the generality in which it is defined.

**Lemma 3.5:** The map $I \otimes R \xrightarrow{1 \otimes d} I \otimes_R \Omega_{R/I}$ induces a surjection $HH_1(R, I) \to I \otimes_R \Omega_{R/I}$.

**Proof:** Define a map from $I \otimes R + R \otimes I$ to $I \otimes_R \Omega_{R/I}$ by sending $x \otimes r$ to $x \otimes dr$ and $r \otimes x$ to $-x \otimes dr$ ($x \in I$, $r \in R$). This map is clearly onto. Using the presentation (3.0), we see that it factors through $HH_1(R, I)$. ∎

**Theorem 3.6:** Suppose $f: R \to S$ is a map of commutative algebras which sends an ideal $I$ of $R$ isomorphically onto an ideal of $S$. Then the map $HH_1(S, I) \to I \otimes_S \Omega_{(S/I)/k}$ of 3.5 induces an isomorphism of $R$-modules

$$HH_0(R, S, I) \cong I \otimes_S \Omega_{S/R} \cong (I/I^2) \otimes_{S/I} \Omega_{(S/I)/(R/I)}.$$

**Proof:** Since $R$ and $S$ are commutative, $HH_0(R, I) = HH_0(S, I) = I$. Therefore $HH_0(R, S, I)$ is the cokernel of the map $HH_1(R, I) \to HH_1(S, I)$. The presentation (3.0) of $HH_1(S, I)$ shows that $HH_0(R, S, I)$ is the quotient of $(S \otimes I + I \otimes S)/(R \otimes I + I \otimes R)$ by the boundary of $S \otimes S \otimes I + S \otimes I \otimes S + I \otimes S \otimes S$. Using $b(1 \otimes x \otimes y) = x \otimes y - 1 \otimes xy + y \otimes x$ for $x \in I$ and $y \in S$, we can eliminate the terms coming from $I \otimes S$. The image of $S \otimes I^2$ vanishes in $HH_0(R, S, I)$ because, for $x \in S$ and $y, z \in I$, we have $b(x \otimes y \otimes z) \equiv -x \otimes yz \pmod{R \otimes I}$. Elementary manipulations now show that

$$HH_0(R, S, I) = \frac{I/I^2 \otimes_k S/R}{b(I \otimes S \otimes S)}.$$

From here the result is a straightforward calculation, the details of which are given in [GW1], 4.1.2 and 4.3. ∎

We need one final calculation before we can prove Theorem 1.2. Consider the principal ideal algebra $B = K[s]/(s^n)$, with $K$ a finite separable extension of $k$. Set

$$(3.7.1) \quad \eta = 1 \otimes s^{n-1} \otimes s + s \otimes s^{n-2} \otimes s + \cdots + s^{i-1} \otimes s^{n-i} \otimes s + \cdots + s^{n-2} \otimes s \otimes s.$$

This is an element of $B \otimes B \otimes B$ whose Hochschild boundary is $n \, s^{n-1} \otimes s$. The following result may be proven in many ways, for example by brute force or symbolic manipulation; our proof is by citation.

**Proposition 3.7:** Let $B = K[s]/(s^n)$, where $K$ is a finite separable extension of $k$. Then

$$HH_2(B) = \begin{cases} K[s]/(s^{n-1}) & \text{on generator } t = s\eta \quad \text{if } \frac{1}{n} \in k \\ K[s]/(s^n) & \text{on generator } t = \eta \quad \text{if } n = 0 \text{ in } k \end{cases}$$

**Proof:** By the Künneth formula ([W, 9.4.1]), there is no loss of generality in assuming that $K = k$. If $\text{char}(k) = 0$, this is exactly the calculation of [GRW, 1.10]. If $\text{char}(k) = p > 0$, the method of *loc. cit.* carries over to yield the result cited. The only subtle point is that if we consider $B$ as a DGA with $s$ in degree 2, then the Eilenberg-Moore spectral sequence degenerates to yield the calculation of [GRW, 2.3]:

$$HH^*_{DG}(B; K) \cong H^*(\Omega \mathbf{CP}^{n-1}; K) \cong K[u, t].$$

These isomorphisms are independent of the characteristic of $K$ because $\mathbf{CP}^{n-1}$ is a formal space ([W, 9.9.12]). ∎

**Proof of Theorem 1.2:** (Cf. [Wa, p.193].) Since $S$ is locally a principal ideal ring, we can choose local generators $s_i \in S$ for the primes over $I$. Then there are integers $n_i$ so that $I$ is locally generated by $\prod s_i^{n_i}$ and $S/I$ is a finite product of the truncated polynomial rings $B_i = K_i[\bar{s}_i]/(\bar{s}_i^{n_i})$. Form the elements $\eta_i = 1 \otimes s_i^{n_i-1} \otimes s_i + \cdots$ of $S \otimes S \otimes S$ corresponding to the elements $\bar{\eta}_i$ of $B_i \otimes B_i \otimes B_i$ described in 3.7.1. Applying the Hochschild boundary in $C_*(S)$ to $s_i^m \eta_i$, we see that

$$b(s_i^m \eta_i) = n_i(s_i^{m+n_i-1} \otimes s_i) - s_i^m \otimes s_i^{n_i}.$$

This element lies in $I \otimes S + S \otimes I$ if $m \geq 1$ (and $m = 0$ if $n_i = 0$ in $k$), and represents the image of $s_i^m \bar{\eta}_i$ under $HH_2(S/I) \to HH_1(S,I)$. Passing to $I \otimes_S \Omega_{S/I}$ as in Lemma 3.5, we obtain the elements

$$(m + n_i) s_i^{m+n_i-1} \otimes ds_i \ \in I \otimes_S \Omega_{S/I} = (I/I^2) \otimes_{S/I} \Omega_{S/I}.$$

Because of our choices of $s_i$, the $B_i$-component of $I \otimes_S \Omega_{S/I}$ is generated by $s_i^{n_i} \otimes ds_i$. Thus a $K_i$-basis of this component is the set of all $s_i^{m+n_i} \otimes ds_i$ with $0 \leq m \leq n_i - 2$ (and $m = n_i - 1$ if $n_i = 0$ in $k$). When $\operatorname{char}(k) = 0$, these are all in the image of $HH_2(S/I)$. Now apply Theroem 3.6. ∎

**Porism 3.8:** When $k$ is a perfect field of characteristic $p \neq 0$, the proof shows that the cokernel of $HH_2(S/I) \to I \otimes_S \Omega_{S/I}$ is the sum of all terms $K_i(s_i^{a-1} \otimes ds_i)$ so that $a \equiv 0 \pmod{p}$ and $n_i + 1 \leq a \leq 2n_i - 1$ (and the term $K_i(s_i^{2n_i-1} \otimes ds_i)$ if $n_i = 0$ in $k$). Thus the map $HH_2(S/I) \to I \otimes_S \Omega_{(S/I)/k}$ is onto if and only if $p > 2n_i$ for all $i$. Since $\Omega_{S/R}$ is a quotient of $\Omega_{(S/I)/k}$, having $p > 2n_i$ is a sufficient, but not necessary, condition for the conclusion of 1.2 to hold.

### Acknowledgements

The authors are grateful to W. Vasconcelos for drawing their attention to Berger's Conjecture, and to B. Osofsky for helpful discussions. The second author also wishes to thank L. Roberts and Queen's University for their kind hospitality during both the research and writing phases of this paper, and to acknowledge the support of NSERC by providing travel funds through a grant to L. Roberts.